\def\titlep{C$^{*}$-bialgebra defined by the
direct sum of Cuntz-Krieger algebras}
\documentclass[11pt]{article}
\usepackage{graphicx,ifthen}
\usepackage{amssymb}

\font\germ=eufm10 at12pt

\def\goth#1{\hbox{\germ#1}}

\newcommand{\qed}{\hbox{\rule[-2pt]{3pt}{6pt}}}
\newcommand{\qedh}{\hfill\qed \\}

\newcommand{\vep}{\varepsilon}

\newtheorem{Thm}{Theorem}[section]

\newtheorem{ex}[Thm]{Example}
\newtheorem{defi}[Thm]{Definition}
\newtheorem{lem}[Thm]{Lemma}

\newtheorem{prop}[Thm]{Proposition}

\def\cal#1{\mathcal #1}
\def\con{{\cal O}_{n}}
\def\coni{{\cal O}_{\infty}}

\def\nset#1{\{1,\ldots,n\}^{#1}}

\def\coa{{\cal O}_{A}}
\def\co#1{{\cal O}_{#1}}
\def\ck{Cuntz-Krieger}

\def\brl{branching law}
\def\bfsnl{{\rm BFS}_{N}(\Lambda)}

\addtocounter{footnote}{1}
\def\sftt#1{
\setcounter{equation}{0}
\addtocounter{footnote}{1}
\section{#1}
}

\def\ssft#1{\subsection{#1}}

\def\sssft#1{\subsubsection{#1}}

\begin{document}
%
%
\def\autherp{Katsunori Kawamura}
\def\emailp{e-mail: kawamura@kurims.kyoto-u.ac.jp.}
\def\addressp{{\small {\it College of Science and Engineering Ritsumeikan University,}}\\
{\small {\it 1-1-1 Noji Higashi, Kusatsu, Shiga 525-8577, Japan}}
}

\def\infw{\Lambda^{\frac{\infty}{2}}V}
\def\zhalfs{{\bf Z}+\frac{1}{2}}
\def\ems{\emptyset}
\def\pmvac{|{\rm vac}\!\!>\!\! _{\pm}}
\def\vac{|{\rm vac}\rangle _{+}}
\def\dvac{|{\rm vac}\rangle _{-}}
\def\ovac{|0\rangle}
\def\tovac{|\tilde{0}\rangle}
\def\expt#1{\langle #1\rangle}
\def\zph{{\bf Z}_{+/2}}
\def\zmh{{\bf Z}_{-/2}}
\def\brl{branching law}
\def\bfsnl{{\rm BFS}_{N}(\Lambda)}
\def\scm#1{S({\bf C}^{N})^{\otimes #1}}
\def\mqb{\{(M_{i},q_{i},B_{i})\}_{i=1}^{N}}
\def\zhalf{\mbox{${\bf Z}+\frac{1}{2}$}}
\def\zmha{\mbox{${\bf Z}_{\leq 0}-\frac{1}{2}$}}
\newcommand{\mline}{
\noindent
\thicklines
\setlength{\unitlength}{.1mm}
\begin{picture}(1000,5)
\put(0,0){\line(1,0){1250}}
\end{picture}
\par
 }
\def\ptimes{\otimes_{\varphi}}
\def\delp{\Delta_{\varphi}}
\def\delps{\Delta_{\varphi^{*}}}
\def\gamp{\Gamma_{\varphi}}
\def\gamps{\Gamma_{\varphi^{*}}}
\def\sem{{\sf M}}
\def\hdelp{\hat{\Delta}_{\varphi}}
\def\tilco#1{\tilde{\co{#1}}}
\def\ndm#1{{\bf M}_{#1}(\{0,1\})}
\def\tndm#1{\tilde{{\bf M}}_{#1}(\{0,1\})}
\def\sck{{\sf CK}_{*}}
\def\hdel{\hat{\Delta}}
\def\ba{\mbox{\boldmath$a$}}
\def\bb{\mbox{\boldmath$b$}}
\def\bc{\mbox{\boldmath$c$}}
\def\be{\mbox{\boldmath$e$}}
\def\bu{\mbox{\boldmath$u$}}
\def\bv{\mbox{\boldmath$v$}}
\def\bw{\mbox{\boldmath$w$}}
\def\bx{\mbox{\boldmath$x$}}
\def\by{\mbox{\boldmath$y$}}
\def\N{{\bf N}}
%
%
%
\setcounter{section}{0}
\setcounter{footnote}{0}
\setcounter{page}{1}
\pagestyle{plain}

%
%
\title{\titlep}
\author{\autherp\thanks{\emailp}
\\
\addressp}
\maketitle

%
%
\begin{abstract}
Let ${\sf CK}_{*}$ denote the C$^{*}$-algebra 
defined by the direct sum of all Cuntz-Krieger algebras.
We introduce a comultiplication $\Delta_{\varphi}$ and a counit 
$\varepsilon$ on ${\sf CK}_{*}$
such that $\Delta_{\varphi}$ is a nondegenerate $*$-homomorphism from 
${\sf CK}_{*}$ to ${\sf CK}_{*}\otimes {\sf CK}_{*}$
and $\varepsilon$ is a $*$-homomorphism from ${\sf CK}_{*}$ to ${\bf C}$.
From this, ${\sf CK}_{*}$ is a counital non-commutative non-cocommutative 
C$^{*}$-bialgebra.
Furthermore, C$^{*}$-bialgebra automorphisms, 
a tensor product of representations
and C$^{*}$-subbialgebras of ${\sf CK}_{*}$ are investigated.
\end{abstract}


%
%
\sftt{Introduction}
\label{section:first}
We have studied C$^{*}$-bialgebras and their construction method.
In this paper, we construct a concrete C$^{*}$-bialgebra
by using well-known C$^{*}$-algebras
and $*$-homomorphisms among them.
We start with our motivation.
%
%
\ssft{Motivation}
\label{subsection:firstone}
In \cite{TS02}, we constructed the C$^{*}$-bialgebra $\co{*}$
defined by the direct sum of all Cuntz algebras except $\co{\infty}$:
%
%
\begin{equation}
\label{eqn:cuntzone}
\co{*}=\co{1}\oplus\co{2}\oplus\co{3}\oplus\co{4}\oplus\cdots
\end{equation}
where $\co{1}$ denotes the $1$-dimensional C$^{*}$-algebra ${\bf C}$
for convenience.
The C$^{*}$-bialgebra $\co{*}$ is non-commutative and non-cocommutative.
We investigated a Haar state, KMS states, C$^{*}$-bialgebra automorphisms,
C$^{*}$-subbialgebras and a comodule-C$^{*}$-algebra of $\co{*}$.
This study was motivated by a certain tensor product 
of representations of Cuntz algebras \cite{TS01}.

Furthermore, we showed a general method to construct
such a C$^{*}$-bialgebra from
a system of C$^{*}$-algebras and $*$-homomorphisms among them.
This system is given as a set
$\{({\cal A}_{a},\varphi_{a,b}):a,b\in \sem\}$
where ${\cal A}_{a}$ is a unital C$^{*}$-algebra,
$\varphi_{a,b}\in {\rm Hom}({\cal A}_{a},{\cal A}_{b})$
and $\sem$ is a monoid ($=$ a semigroup with unit), and 
they satisfy several assumptions.
The system will be explained in $\S$ \ref{section:third}.
For the case of $\co{*}$,
$\sem$ is the monoid ${\bf N}=\{1,2,3,\ldots\}$
with respect to the multiplication.
Although ${\bf N}$ is commutative,
the commutativity of $\sem$ is not necessary for the system 
$\{({\cal A}_{a},\varphi_{a,b}):a,b\in \sem\}$ in general. 
Hence we are interesting in an example of C$^{*}$-bialgebra associated with
such a system over a non-commutative monoid $\sem$.

On the other hand,
Cuntz-Krieger algebras are well-known natural generalizations of 
Cuntz algebras, and
many studies about Cuntz algebras are
generalized to Cuntz-Krieger algebras.
Therefore the idea of a C$^{*}$-bialgebra associated with
Cuntz-Krieger algebras like $\co{*}$ in (\ref{eqn:cuntzone}) is natural.

In this paper,
we construct a C$^{*}$-bialgebra 
by using Cuntz-Krieger algebras as an example
of system over a non-commutative monoid
consisting of all nondegenerate matrices with entries $0$ or $1$.
The monoid structure will be explained in $\S$ \ref{subsection:firstthree}
and $\S$ \ref{section:second}.

%
%
\ssft{C$^{*}$-bialgebras}
\label{subsection:firsttwo}
We prepare terminology about C$^{*}$-bialgebra according to \cite{ES,KV,MNW}.
Since we use many algebras and matrices at once in this paper, 
we write ${\cal A},{\cal B},{\cal C},\ldots$ as (C$^{*}$-, co-, bi-) algebras 
and write $A,B,C.\ldots$ as matrices.
For two C$^{*}$-algebras ${\cal A}$ and ${\cal B}$,
we write ${\rm Hom}({\cal A},{\cal B})$ as the set of all $*$-homomorphisms 
from ${\cal A}$ to ${\cal B}$.
Assume that every tensor product $\otimes$ as below means the minimal C$^{*}$-tensor product.
%
%
\begin{defi}
\label{defi:cstar}
A pair $({\cal A},\Delta)$ is a C$^{*}$-bialgebra
if ${\cal A}$ is a C$^{*}$-algebra and $\Delta\in {\rm Hom}({\cal A},M({\cal A}\otimes {\cal A}))$ 
where $M({\cal A}\otimes {\cal A})$ denotes the multiplier algebra of ${\cal A}\otimes {\cal A}$
such that the linear span of $\{\Delta(a)(b\otimes c):a,b,c\in {\cal A}\}$ 
is norm dense in ${\cal A}\otimes {\cal A}$ and 
the following holds:
%
%
\begin{equation}
\label{eqn:bialgebratwo}
(\Delta\otimes id)\circ \Delta=(id\otimes\Delta)\circ \Delta.
\end{equation}
We call $\Delta$ the comultiplication of ${\cal A}$.
\end{defi}

\noindent
We state that a C$^{*}$-bialgebra $({\cal A},\Delta)$ is {\it strictly proper} 
if $\Delta(a)\in {\cal A}\otimes {\cal A}$ for any $a\in {\cal A}$;
$({\cal A},\Delta)$ is {\it unital}
if ${\cal A}$ is unital and $\Delta$ is unital;
$({\cal A},\Delta)$ is {\it counital}
if there exists $\vep\in {\rm Hom}({\cal A},{\bf C})$ such that
%
%
\begin{equation}
\label{eqn:counit}
(\vep\otimes id)\circ \Delta\cong id \cong (id\otimes \vep)\circ\Delta.
\end{equation}
We call $\vep$ the {\it counit} of ${\cal A}$ and write $({\cal A},\Delta,\vep)$ 
as the counital C$^{*}$-bialgebra $({\cal A},\Delta)$ with the counit $\vep$.
Remark that Definition \ref{defi:cstar} does not mean 
$\Delta({\cal A})\subset {\cal A}\otimes {\cal A}$.
If ${\cal A}$ is unital, then  $({\cal A},\Delta)$ is strictly proper.
A {\it bialgebra} in the purely algebraic theory \cite{Abe,Kassel} means 
a unital counital strictly proper bialgebra with the unital counit
with respect to the algebraic tensor product,
which does not need to have an involution.
Hence C$^{*}$-bialgebra is not a bialgebra in general.

%
%
\ssft{Main theorems}
\label{subsection:firstthree}
In this subsection, we state our main theorems.
A matrix $A$ is {\it nondegenerate} if any column and any row are not zero.
For $1\leq n<\infty$,
let $\ndm{n}$ denote the set of all
nondegenerate $n\times n$ matrices with entries $0$ or $1$.
In particular, $\ndm{1}=\{1\}$.
Define
%
%
\begin{equation}
\label{eqn:nondegenerate}
\ndm{*}\equiv \cup \{\ndm{n}:n\in {\bf N}\}.
\end{equation}
For $A\in\ndm{n}$,
let  $\coa$ denote the Cuntz-Krieger algebra by $A$ with
the canonical generators $s_{1}^{(A)},\ldots,s_{n}^{(A)}$ 
where we define $s_{1}^{(1)}=I_{1}$ 
and $I_{1}$ is the unit of the C$^{*}$-algebra $\co{1}={\bf C}$.

For $A=(a_{ij})\in \ndm{n}$ and $B=(b_{ij})\in\ndm{m}$,
define the Kronecker product $A\boxtimes B\in\ndm{nm}$ 
of $A$ and $B$ by
%
%
\begin{equation}
\label{eqn:phitwo}
(A\boxtimes B)_{m(i-1)+j,m(i^{'}-1)+j^{'}}\equiv a_{ii^{'}}b_{jj^{'}}
\end{equation}
for $i,i^{'}\in \nset{}$ and $j,j^{'}\in\{1,\ldots,m\}$ \cite{Dief}.
In addition,
define the map $\varphi_{A,B}$ from $\co{A\boxtimes B}$ 
to the minimal tensor product $\co{A}\otimes \co{B}$ by
%
%
\begin{equation}
\label{eqn:varphi}
\varphi_{A,B}(s_{m(i-1)+j}^{(A\boxtimes B)})
\equiv s_{i}^{(A)}\otimes s_{j}^{(B)}\quad 
(i\in\nset{},\,j\in\{1,\ldots,m\}).
\end{equation}
Then we can verify that $\varphi_{A,B}$ 
is uniquely extended to a unital $*$-embedding of 
$\co{A\boxtimes B}$ into $\co{A}\otimes \co{B}$.
Remark that 
$\co{A\boxtimes B}\not\cong\co{A}\otimes \co{B}$
in general.
%
%
\begin{Thm}
\label{Thm:maintsfour}
Define the C$^{*}$-algebra $\sck$ by the direct sum of $\{\coa:A\in\ndm{*}\}$:
%
%
\begin{equation}
\label{eqn:main}
\sck\equiv \oplus\{\coa:A\in \ndm{*}\}.
\end{equation}
Let $\sck\otimes \sck$ denote the minimal tensor product of $\sck$.
For the set $\varphi=\{\varphi_{A,B}:A,B\in \ndm{*}\}$
in (\ref{eqn:varphi}),
define $\delp\in {\rm Hom}(\sck,\sck\otimes \sck)$ and
$\vep\in{\rm Hom}({\sf CK}_{*},{\bf C})$ by
\[
\delp \equiv \oplus\{\delp^{(A)}:A\in \ndm{*}\},\quad
\delp^{(A)}(x)\equiv 
\sum_{(B,C)\in {\cal N}_{A}}\varphi_{B,C}(x)\quad(x\in \coa),
\]
\[\vep(x)\equiv 0\quad(x\in \oplus\{\coa:A\in \ndm{*},\,A\ne 1\}),\quad
\vep(x)\equiv x\quad(x\in \co{1})\]
where
%
%
\begin{equation}
\label{eqn:na}
{\cal N}_{A}\equiv \{(B,C)\in \ndm{*}\times \ndm{*}:B\boxtimes C=A\}.
\end{equation}
Then $({\sf CK}_{*},\delp,\vep)$ is a counital C$^{*}$-bialgebra.
\end{Thm}

\noindent
By construction,
$({\sf CK}_{*},\delp,\vep)$ is strictly proper and nonunital.
The essential part of Theorem \ref{Thm:maintsfour} 
for the bialgebra structure
is the set $\varphi=\{\varphi_{A,B}:A,B\in \ndm{*}\}$
in (\ref{eqn:varphi}).
The nontrivial fact is that
the set of all Cuntz-Krieger algebras
has such a set of $*$-embeddings.

We show several properties of $({\sf CK}_{*},\delp,\vep)$ as follows.
%
%
\begin{Thm}
\label{Thm:maintwo}
\begin{enumerate}
\item
Let $\widetilde{{\sf CK}}_{*}$ denote the smallest unitization of ${\sf CK}_{*}$
as a C$^{*}$-algebra.
Then there exists a comultiplication $\hdelp$ and 
a counit $\tilde{\vep}$ of $\widetilde{{\sf CK}}_{*}$ such that
$(\widetilde{{\sf CK}}_{*},\hdelp,\tilde{\vep})$
is a strictly proper unital counital C$^{*}$-bialgebra,
$\hdelp|_{{\sf CK}_{*}}=\delp$ and $\tilde{\vep}|_{{\sf CK}_{*}}=\vep$.
\item
There exists a dense (unital counital $*$-) subbialgebra ${\goth A}_{0}$ 
of $\widetilde{{\sf CK}}_{*}$ such that $\hdelp({\goth A}_{0})$
is included in the algebraic tensor product ${\goth A}_{0}\odot{\goth A}_{0}$
of ${\goth A}_{0}$.
\item
There is no antipode for any dense subbialgebra of ${\sf CK}_{*}$.
\item
Define $F_{n}\in\ndm{n}$ by $(F_{n})_{ij}=1$ for any $i,j$.
Then 
\[{\sf C}_{*}\equiv  \oplus\{\co{F_{n}}:n\in {\bf N}\}\]
is a counital C$^{*}$-subbialgebra of ${\sf CK}_{*}$
which is isomorphic to $\co{*}$ in (\ref{eqn:cuntzone})
as a C$^{*}$-bialgebra.
\end{enumerate}
\end{Thm}

We discuss results in Theorem \ref{Thm:maintsfour}
and \ref{Thm:maintwo} here.
\begin{enumerate}
\item
The bialgebra structure does not appear unless one takes
the direct sum of $\coa$'s.
It is a rare example that every Cuntz-Krieger algebras appear all at once.
\item
The bialgebra $({\goth A}_{0},\hat{\Delta}_{\varphi}|_{{\goth A}_{0}},
\tilde{\vep}|_{{\goth A}_{0}})$ 
in Theorem \ref{Thm:maintwo} (ii) is constructed by using 
neither C$^{*}$-algebra nor the C$^{*}$-tensor product.
(This will be shown in the proof of Theorem \ref{Thm:maintwo} 
in $\S$ \ref{section:third}.)
Furthermore,
this kind of bialgebra has not been known 
in the purely algebraic theory of bialgebras yet.
Hence 
$({\goth A}_{0},\hat{\Delta}_{\varphi}|_{{\goth A}_{0}},
\tilde{\vep}|_{{\goth A}_{0}})$ is a remarkable example
not only in operator algebras but also in the purely algebraic theory. 
\item
Since there is no standard comultiplication of ${\sf CK}_{*}$,
$({\sf CK}_{*},\delp)$ is not a deformation 
of any cocommutative C$^{*}$-bialgebra. 
\item
It is clear that ${\sf C}_{*}$ is isomorphic to
$\co{*}$ as a ``C$^{*}$-algebra" in Theorem \ref{Thm:maintwo} (iv),
but it is not trivial whether they are isomorphic 
or not as a ``C$^{*}$-bialgebra."
Since $\co{*}$ is non-cocommutative,
so is $({\sf CK}_{*},\delp,\vep)$.
\end{enumerate}

\noindent
{\bf Question}
In \cite{TS02}, we show that 
$\coni$ is a comodule-C$^{*}$-algebra of the C$^{*}$-bialgebra $\co{*}$ 
in (\ref{eqn:cuntzone}).
Find a comodule-C$^{*}$-algebra of ${\sf CK}_{*}$.
We guess that Cuntz-Krieger algebras for
infinite matrices \cite{EL,Fang} are the candidates.\\

In $\S$ 2, we explain properties of the monoid of all nondegenerate matrices.
In $\S$ 3, we show a general method to construct C$^{*}$-bialgebras.
Proofs of Theorem \ref{Thm:maintsfour} and \ref{Thm:maintwo} 
will be given in $\S$ 3.4.
In $\S$ 4, C$^{*}$-bialgebra automorphisms of ${\sf CK}_{*}$
and a tensor product of representations of Cuntz-Krieger algebras are discussed.
In $\S$ 5, we show examples of C$^{*}$-subbialgebras of ${\sf CK}_{*}$.

%
%
\sftt{Monoid of matrices}
\label{section:second}
In this section, we prepare a monoid of all nondegenerate matrices
with respect to the Kronecker product
in order to construct the C$^{*}$-bialgebra in Theorem \ref{Thm:maintsfour}.
A {\it monoid} is a set $\sem$ equipped with a binary associative operation, 
$\sem\times \sem\ni(a,b)\mapsto ab\in \sem$,
and a unit with respect to the operation.

%
%
\ssft{Kronecker product of vectors}
\label{subsection:secondone}
We introduce a realization of the tensor product of vectors.
Let $\{e_{i}^{(n)}\}_{i=1}^{n}$ denote
the standard basis of the complex vector space ${\bf C}^{n}$ for $1\leq n<\infty$.
Define the linear isomorphism $\Psi_{n,m}$ 
from ${\bf C}^{n}\otimes {\bf C}^{m}$ to ${\bf C}^{nm}$ by
\[\Psi_{n,m}(e_{i}^{(n)}\otimes e_{j}^{(m)})\equiv e^{(nm)}_{m(i-1)+j}
\quad((i,j)\in \{1,\ldots,n\}\times \{1,\ldots,m\}).\]
We write
%
%
\begin{equation}
\label{eqn:vector}
v\boxtimes w\equiv \Psi_{n,m}(v\otimes w)\quad
(v\in {\bf C}^{n},\,w\in {\bf C}^{m}).
\end{equation}
We call $\boxtimes$ the {\it Kronecker product}
of $\{{\bf C}^{n}:n\geq 1\}$.
For any $n,m,l\geq 1$, the following holds in ${\bf C}^{nml}$:
\[(v\boxtimes w)\boxtimes u=v\boxtimes (w\boxtimes u)\quad
(v\in {\bf C}^{n},\,w\in {\bf C}^{m},\,u\in {\bf C}^{l}).\]
Since $e_{2}^{(2)}\boxtimes e_{2}^{(3)}=e_{5}^{(6)}$ 
and $e_{2}^{(3)}\boxtimes e_{2}^{(2)}=e_{4}^{(6)}$,
$v\boxtimes w\ne w\boxtimes v$ in general.

Let ${\bf Z}^{n}$ denote the subset of ${\bf C}^{n}$
consisting of $\sum_{i=1}^{n}a_{i}e_{i}^{(n)}$ for
$a_{i}\in{\bf Z}$ for each $i$.
Then we see that the set $\bigcup_{n\geq 1}{\bf Z}^{n}$
is closed with respect to the operation $\boxtimes$. 

%
%
\ssft{Monoid of all matrices}
\label{subsection:secondtwo}
We review properties of the Kronecker product $\boxtimes$ in (\ref{eqn:phitwo}).
Define the set ${\cal M}$ by the set-theoretical union
of $\{M_{n}({\bf C}):n\geq 1\}$:
%
%
\begin{equation}
\label{eqn:monoidone}
{\cal M}\equiv \bigcup_{n\geq 1}M_{n}({\bf C}).
\end{equation}
Then the following holds:
\[
\begin{array}{ll}
A\boxtimes 1=1\boxtimes A=A\quad&(A\in {\cal M}),\\\
\\
(A\boxtimes B)\boxtimes C=A\boxtimes (B\boxtimes C)
\quad&(A,B,C\in {\cal M})\\
\end{array}
\]
where $1\in M_{1}({\bf C})={\bf C}$.
We can regard ${\bf C}^{n}$ as the subset of 
$M_{n}({\bf C})$ which consists of all diagonal matrices.
Then $\boxtimes$ in (\ref{eqn:vector}) coincides with
the Kronecker product $\boxtimes$ of matrices.
From this, $({\cal M},\boxtimes)$ is a non-commutative monoid with the unit $1$.
For $A,B\in M_{n}({\bf C})$ and $C,D\in M_{m}({\bf C})$,
we see that
%
%
\begin{equation}
\label{eqn:multiplication}
(A\boxtimes C)(B\boxtimes D)=AB\boxtimes CD.
\end{equation}
Let $1_{n}$ denote the identity matrix of $M_{n}({\bf C})$.
Then $1_{n}\boxtimes 1_{m}=1_{nm}$.

A matrix $A\in M_{n}({\bf C})$ is {\it irreducible} if 
for any $i,j\in\{1,\ldots,n\}$,
there exists $k\in {\bf N}$ such that
$(A^{k})_{i,j}\ne 0$ where $A^{k}=A\cdots A$ ($k$-times).
From (\ref{eqn:multiplication}),
the following holds.
%
%
\begin{lem}
\label{lem:matrix}
For $A,B\in {\cal M}$,
both $A$ and $B$ are irreducible if and only if $A\boxtimes B$
is irreducible.
\end{lem}

For $v\in {\bf C}^{n}$ and $w\in {\bf C}^{m}$,
let $v\boxtimes w$ be as in (\ref{eqn:vector}).
Then
$(A\boxtimes B)(v\boxtimes w)=Av\boxtimes Bw$
for any $A\in M_{n}({\bf C})$ and $B\in M_{m}({\bf C})$.
From this, the following inclusion holds:
%
%
\begin{equation}
\label{eqn:kernel}
\ker(1_{n}-A)\boxtimes \ker(1_{m}-B)\subset
\ker(1_{nm}-A\boxtimes B).
\end{equation}

%
%
\ssft{Monoid of nondegenerate $0$-$1$ matrices}
\label{subsection:secondthree}
Let $\ndm{*}$ be as in (\ref{eqn:nondegenerate}).
Then we see that $\ndm{*}$ is a non-commutative submonoid 
of $({\cal M},\boxtimes)$ in (\ref{eqn:monoidone}) with the common unit $1\in \ndm{1}$.
For $A,A^{'}\in \ndm{n}$ and $B,B^{'}\in \ndm{m}$,
if $A\boxtimes B=A^{'}\boxtimes B^{'}$, then $A=A^{'}$ and $B=B^{'}$.
Let $F_{n}$ be as in Theorem \ref{Thm:maintwo} (iv).
We see that $F_{n}\boxtimes F_{m}=F_{nm}$. 
From this, a decomposition of an element $A$ of $\ndm{*}$ 
with respect to $\boxtimes$ is not unique in general.
The map $F$ from ${\bf N}$ to $\ndm{*}$ is an embedding 
of the monoid $({\bf N},\cdot)$ into $(\ndm{*},\boxtimes)$. 
If $A\in \ndm{n}$, then $A{\bf Z}^{n}\subset {\bf Z}^{n}$.
From this, we regard $A\in \ndm{n}$ 
as a map from ${\bf Z}^{n}$ to ${\bf Z}^{n}$.
We see that $A\boxtimes B=F_{n}$
if and only if there exist $m,l\in {\bf N}$
such that $ml=n$ and $A=F_{m}$ and $B=F_{l}$.

Let ${\bf IM}$ denote the set of all irreducible matrices in $\ndm{*}$.
From Lemma \ref{lem:matrix},
${\bf IM}$ is a submonoid of $(\ndm{*},\boxtimes)$
such that any divisor of any element in ${\bf IM}$ belongs to ${\bf IM}$.
%
%
\begin{lem}
\label{lem:permutation}
For $A,B\in {\bf IM}$,
if neither $A$ nor $B$ are permutation matrices,
then  $A\boxtimes B$ is not a permutation matrix.
\end{lem}

%
%
\sftt{C$^{*}$-weakly coassociative system}
\label{section:third}
In this section, we review the general method to construct C$^{*}$-bialgebras \cite{TS02},
and prove Theorem \ref{Thm:maintsfour}
and \ref{Thm:maintwo}.

%
%
\ssft{C$^{*}$-algebras}
\label{subsection:thirdone}
In this subsection, we explain  basic facts of
the direct sum of general C$^{*}$-algebras and Cuntz-Krieger algebras.
%
%
\sssft{Direct sum}
\label{subsubsection:thirdoneone}
For an infinite set $\{{\cal A}_{i}:i\in \Omega\}$ of C$^{*}$-algebras,
we define the C$^{*}$-algebra
$\bigoplus_{i\in\Omega} {\cal A}_{i}$ as follows \cite{Blackadar2006}:
\[
\bigoplus_{i\in\Omega} {\cal A}_{i}\equiv 
\{(a_{i}):\|(a_{i})\|\to 0\mbox{ as }i\to\infty\}
\]
in the sense that for every $\vep>0$ there are only finitely many
$i$ for which $\|a_{i}\|>\vep$.
We call $\bigoplus_{i\in\Omega} {\cal A}_{i}$ 
the {\it direct sum} of ${\cal A}_{i}$'s. 
The algebraic direct sum $\oplus_{alg}\{{\cal A}_{i}:i\in\Omega\}$ 
is a dense $*$-subalgebra of $\oplus\{{\cal A}_{i}:i\in\Omega\}$.

Let $\{{\cal B}_{i}:i\in \Omega\}$ be another set of C$^{*}$-algebras
and let $\{f_{i}:i\in \Omega\}$ be a set of $*$-homomorphisms such that
$f_{i}\in {\rm Hom}({\cal A}_{i},{\cal B}_{i})$ for each $i\in \Omega$.
Then we obtain 
$\oplus_{i\in\Omega} f_{i}\in
{\rm Hom}(\oplus_{i\in\Omega} {\cal A}_{i},\oplus_{i\in\Omega} {\cal B}_{i})$.
If $f_{i}$ is nondegenerate for each $i$,
then $\oplus_{i\in \Omega}f_{i}$ is also nondegenerate.
If both ${\cal A}_{i}$ and ${\cal B}_{i}$ 
are unital and $f_{i}$ is unital for each $i\in \Omega$,
then $\oplus_{i\in \Omega}f_{i}$ is nondegenerate.

%
%
\sssft{Cuntz-Krieger algebras}
\label{subsubsection:thirdonetwo}
Let $\ndm{*}$ be as in (\ref{eqn:nondegenerate}).
For $A=(a_{ij})\in \ndm{n}$,
$\coa$ is the
{\it \ck\ algebra by $A$} if 
$\coa$ is a C$^{*}$-algebra
which is universally generated by
partial isometries $s_{1},\ldots,s_{n}$
and they satisfy $s_{i}^{*}s_{i}=\sum_{j=1}^{n}a_{ij}s_{j}s_{j}^{*}$
for $i=1,\ldots,n$ and $\sum_{i=1}^{n}s_{i}s_{i}^{*}=I$ \cite{CK}.
Especially,
when $a_{ij}=1$ for each $i,j=1,\ldots,n$,
$\coa$ is $*$-isomorphic to the {\it Cuntz algebra} $\con$ \cite{C}.
Let $F_{n}$ be as in Theorem \ref{Thm:maintwo} (iv).
Then $\co{F_{n}}\cong \con$. 
The C$^{*}$-algebra $\coa$ is simple 
if and only if $A$ is irreducible and not a permutation matrix.
It is known that the $K_{1}$-group of $\coa$
is isomorphic to $\ker(1-A^{t}:{\bf Z}^{n}\to {\bf Z}^{n})$
where $A^{t}$ denotes the transposed matrix of $A$ \cite{Wegge}.
For embeddings of Cuntz-Krieger algebras, see \cite{CK01}. 

Let $\boxtimes$ be as in (\ref{eqn:phitwo}).
From Lemma \ref{lem:permutation},
if both $\coa$ and $\co{B}$ are simple, then so is $\co{A\boxtimes B}$.
From (\ref{eqn:kernel}) and $(A\boxtimes B)^{t}=A^{t}\boxtimes B^{t}$,
$K_{1}(\coa)\otimes_{{\bf Z}}K_{1}(\co{B})$ 
is embedded into $K_{1}(\co{A\boxtimes B})$.

%
%
\ssft{C$^{*}$-bialgebras continued}
\label{subsection:thirdtwo}
This subsection is the succeeding part of $\S$ \ref{subsection:firsttwo}.  
%
%
\sssft{Bialgebras in the purely algebraic theory}
\label{subsubsection:thirdtwoone}
In order to consider subbialgebras 
and bialgebra morphisms for C$^{*}$-bialgebras,
we start with bialgebras in the purely algebraic theory
according to \cite{Kassel}.
In this subsubsection,
any tensor product $\otimes$ means the algebraic tensor product.
Let $k$ be the ground field with the unit $1$.
A {\it coalgebra} is a triplet $({\cal C},\Delta,\vep)$ where
${\cal C}$ is a vector space and $\Delta:{\cal C}\to {\cal C}\otimes {\cal C}$
and $\vep:{\cal C}\to k$ are linear maps satisfying 
axioms (\ref{eqn:bialgebratwo}) and (\ref{eqn:counit}).
A {\it bialgebra} is a quintuple $({\cal B},m,\eta,\Delta,\vep)$
where $({\cal B},m,\eta)$ is a unital associative algebra and 
$({\cal B},\Delta,\vep)$ is a counital coassociative coalgebra
such that both $\Delta$ and $\vep$
are unital algebra morphisms.
An endomorphism $S$ of ${\cal B}$ is called an {\it antipode}
for $({\cal B},m,\eta,\Delta,\vep)$ if $S$ satisfies
$m\circ (id\otimes S)\circ \Delta=\eta\circ \vep=
m\circ (S\otimes id)\circ \Delta$.
If an antipode exists on ${\cal B}$, then it is unique.
If $({\cal B},m,\eta,\Delta,\vep)$ has the antipode $S$, then
$({\cal B},m,\eta,\Delta,\vep,S)$ is called a {\it Hopf algebra}.

For two bialgebras ${\cal A}$ and ${\cal B}$, 
a map $f$ from ${\cal A}$ to ${\cal B}$ is a {\it bialgebra morphism} if 
$f$ is a unital algebra morphism and 
$\Delta_{{\cal B}}\circ f= (f\otimes f)\circ \Delta_{{\cal A}}$ and 
$\vep_{{\cal B}}\circ f=\vep_{{\cal A}}$.
A map $f$ is a {\it bialgebra automorphism} of ${\cal A}$ if $f$ 
is a bijective bialgebra morphism from ${\cal A}$ to ${\cal A}$.

Assume that ${\cal A}$ is an algebra without unit.
Let $\tilde{{\cal A}}$ denote the direct sum $k\oplus {\cal A}$ of  
two $k$-vector spaces $k$ and ${\cal A}$.
By the operation
$(a,x)(b,y)=(ab,ay+bx+xy)$ for $(a,x),(b,y)\in \tilde{{\cal A}}$,
$\tilde{{\cal A}}$ is an algebra with the unit $I_{\tilde{{\cal A}}}=(1,0)$.
By the natural embedding $\iota$ of ${\cal A}$ into $\tilde{{\cal A}}$,
${\cal A}$ is a two-sided ideal of $\tilde{{\cal A}}$ such that
$\tilde{{\cal A}}/{\cal A}\cong k$.
We call $\tilde{{\cal A}}$ the {\it unitization} of ${\cal A}$.
From the proof of Lemma 2.2 in \cite{TS02},
the following is verified.
%
%
\begin{lem}
\label{lem:extensionunittwo}
Assume that ${\cal A}$ is an algebra without unit,
$\Delta\in {\rm Hom}({\cal A},{\cal A}\otimes {\cal A})$
and $\vep\in {\rm Hom}({\cal A},k)$ which satisfy
(\ref{eqn:bialgebratwo}) and (\ref{eqn:counit}).
Then there exist $\hat{\Delta}
\in {\rm Hom}(\tilde{{\cal A}},\tilde{{\cal A}}\otimes \tilde{{\cal A}})$ 
and $\tilde{\vep}\in{\rm Hom}(\tilde{{\cal A}},k)$ such that 
$(\tilde{{\cal A}},\hat{\Delta},\tilde{\vep})$ is a bialgebra,
$\hat{\Delta}|_{{\cal A}}=\Delta$ and $\tilde{\vep}|_{{\cal A}}=\vep$.
\end{lem}

A {\it $*$-bialgebra} over ${\bf C}$ 
is a bialgebra $({\cal B},m,\eta,\Delta,\vep)$ with an involution $*$
such that $({\cal B},m,\eta)$ is a $*$-algebra and both $\Delta$ and $\vep$
are $*$-algebra morphisms \cite{Schurmann}.
For two $*$-bialgebras ${\cal A}$ and ${\cal B}$, 
a map $f$ from ${\cal A}$ to ${\cal B}$ is a {\it $*$-bialgebra morphism} 
if $f$ is a bialgebra morphism and it preserves $*$.

%
%
\sssft{C$^{*}$-bialgebras}
\label{subsubsection:thirdtwotwo}
In this subsubsection,
we assume that $\otimes$ means the minimal tensor product of C$^{*}$-algebras.
For two C$^{*}$-bialgebras $({\cal A}_{1},\Delta_{1})$ 
and $({\cal A}_{2},\Delta_{2})$,
$f$ is a {\it C$^{*}$-bialgebra morphism} from 
$({\cal A}_{1},\Delta_{1})$ to $({\cal A}_{2},\Delta_{2})$ 
if $f$ is a nondegenerate $*$-homomorphism 
from ${\cal A}_{1}$ to $M({\cal A}_{2})$
such that $(f\otimes f)\circ \Delta_{1}=\Delta_{2}\circ f$.
In addition, if $f({\cal A}_{1})\subset {\cal A}_{2}$, then $f$ 
is called {\it strictly proper}.
For two counital C$^{*}$-bialgebras 
$({\cal A}_{1},\Delta_{1},\vep_{1})$ and $({\cal A}_{2},\Delta_{2},\vep_{2})$,
a C$^{*}$-bialgebra morphism $f$ from 
$({\cal A}_{1},\Delta_{1},\vep_{1})$ to $({\cal A}_{2},\Delta_{2},\vep_{2})$ 
is {\it counital} if $\vep_{2}\circ f=\vep_{1}$.
A map $f$ is a {\it C$^{*}$-bialgebra automorphism} of
a C$^{*}$-bialgebra $({\cal A},\Delta)$ if 
$f$ is a bijective C$^{*}$-bialgebra morphism from ${\cal A}$ to ${\cal A}$.
Remark that any $*$-homomorphism among C$^{*}$-algebras
is a morphism of C$^{*}$-algebras, 
but a C$^{*}$-bialgebra morphism is not always a $*$-homomorphism.

A C$^{*}$-algebra ${\cal A}_{0}$ is a {\it C$^{*}$-subbialgebra} of 
a C$^{*}$-bialgebra $({\cal A},\Delta)$ if 
${\cal A}_{0}$ is a C$^{*}$-subalgebra of ${\cal A}$ such that 
the inclusion map is nondegenerate 
and $\Delta({\cal A}_{0})\subset M({\cal A}_{0}\otimes {\cal A}_{0})$.
A C$^{*}$-subbialgebra ${\cal A}_{0}$ 
of a counital C$^{*}$-bialgebra $({\cal A},\Delta,\vep)$ is {\it counital} 
if $({\cal A}_{0},\Delta|_{{\cal A}_{0}},\vep|_{{\cal A}_{0}})$
is a counital C$^{*}$-bialgebra.

%
%
\begin{lem}
\label{lem:extensionunit}
(Lemma 2.2 in \cite{TS02}).
Assume that $({\cal A},\Delta,\vep)$ is a strictly proper nonunital counital
C$^{*}$-bialgebra and $\tilde{{\cal A}}$ 
is the smallest unitization of ${\cal A}$.
Then there exist $\hat{\Delta}
\in {\rm Hom}(\tilde{{\cal A}},\tilde{{\cal A}}\otimes \tilde{{\cal A}})$ 
and $\tilde{\vep}\in{\rm Hom}(\tilde{{\cal A}},{\bf C})$ such that 
$(\tilde{{\cal A}},\hat{\Delta},\tilde{\vep})$ is a 
strictly proper unital counital C$^{*}$-bialgebra
with the unital counit $\tilde{\vep}$,
$\hat{\Delta}|_{{\cal A}}=\Delta$ and $\tilde{\vep}|_{{\cal A}}=\vep$.
\end{lem}

For two counital C$^{*}$-bialgebras 
$({\cal A},\Delta_{{\cal A}},\vep_{{\cal A}})$ and 
$({\cal B},\Delta_{{\cal B}},\vep_{{\cal B}})$,
assume that they satisfy the assumption in Lemma \ref{lem:extensionunit}.
If $f$ is a strictly proper counital C$^{*}$-bialgebra morphism 
from ${\cal A}$ to ${\cal B}$,
then we can verify that
the unique extension $\tilde{f}$ of $f$ 
from $\tilde{{\cal A}}$ to $\tilde{{\cal B}}$
is a unital counital C$^{*}$-bialgebra 
morphism from 
$(\tilde{{\cal A}},\hat{\Delta}_{{\cal A}},\tilde{\vep}_{{\cal A}})$
to $(\tilde{{\cal B}},\hat{\Delta}_{{\cal B}},\tilde{\vep}_{{\cal B}})$.
In particular,
if ${\cal A}$ is a counital C$^{*}$-subbialgebra 
of $({\cal B},\Delta_{{\cal B}},\vep_{{\cal B}})$,
then the inclusion map $\iota$ is extended to
the unital inclusion map $\tilde{\iota}$ of $\tilde{{\cal A}}$ 
into $\tilde{{\cal B}}$.
Hence $\tilde{{\cal A}}$ is a counital unital C$^{*}$-subbialgebra 
of $\tilde{{\cal B}}$.

%
%
\ssft{C$^{*}$-weakly coassociative system}
\label{subsection:thirdthree}
According to $\S$ 3 in \cite{TS02},
we review a general method to construct a C$^{*}$-bialgebra
from a set of C$^{*}$-algebras and $*$-homomorphisms among them.
We assume that $\otimes$ means the minimal tensor product of C$^{*}$-algebras
but not the algebraic tensor product.
%
%
\begin{defi}
\label{defi:axiom}
Let $\sem$ be a monoid with the unit $e$.
A data $\{({\cal A}_{a},\varphi_{a,b}):a,b\in \sem\}$
is a C$^{*}$-weakly coassociative system (= C$^{*}$-WCS) over $\sem$ if 
${\cal A}_{a}$ is a unital C$^{*}$-algebra for $a\in \sem$
and $\varphi_{a,b}$ is a unital $*$-homomorphism
from ${\cal A}_{ab}$ to ${\cal A}_{a}\otimes {\cal A}_{b}$
for $a,b\in \sem$ such that
\begin{enumerate}
\item
for all $a,b,c\in \sem$, the following holds:
%
%
\begin{equation}
\label{eqn:wcs}
(id_{a}\otimes \varphi_{b,c})\circ \varphi_{a,bc}
=(\varphi_{a,b}\otimes id_{c})\circ \varphi_{ab,c}
\end{equation}
where $id_{x}$ denotes the identity map on ${\cal A}_{x}$ for $x=a,c$,
\item
there exists a counit $\vep_{e}$ of ${\cal A}_{e}$ 
such that $({\cal A}_{e},\varphi_{e,e},\vep_{e})$ 
is a counital C$^{*}$-bialgebra,
\item
$\varphi_{e,a}(x)=I_{e}\otimes x$ and
$\varphi_{a,e}(x)=x\otimes I_{e}$ for $x\in {\cal A}_{a}$ and $a\in \sem$.
\end{enumerate}
\end{defi}

\noindent 
From this definition, the following holds.
%
%
\begin{Thm}
\label{Thm:mainthree}
(Theorem 3.1 in \cite{TS02}).
Let $\{({\cal A}_{a},\varphi_{a,b}):a,b\in \sem\}$ be a C$^{*}$-WCS 
over a monoid $\sem$.
Assume that $\sem$ satisfies that 
%
%
\begin{equation}
\label{eqn:finiteness}
\#{\cal N}_{a}<\infty \mbox{ for each }a\in \sem
\end{equation}
where ${\cal N}_{a}\equiv\{(b,c)\in \sem\times \sem:\,bc=a\}$.
Define C$^{*}$-algebras 
\[{\cal A}_{*}\equiv  \oplus \{{\cal A}_{a}:a\in \sem\},\quad
{\cal C}_{a}\equiv 
\oplus \{{\cal A}_{b}\otimes {\cal A}_{c}:(b,c)\in {\cal N}_{a}\}
\quad (a\in\sem). \]
Define $\Delta^{(a)}_{\varphi}\in{\rm Hom}({\cal A}_{a},{\cal C}_{a})$,
$\Delta_{\varphi}
\in {\rm Hom}({\cal A}_{*}, {\cal A}_{*}\otimes {\cal A}_{*})$ and
$\vep\in {\rm Hom}({\cal A}_{*},{\bf C})$ by 
\[\Delta^{(a)}_{\varphi}(x)\equiv \sum_{(b,c)\in {\cal N}_{a}}
\varphi_{b,c}(x)\quad(x\in {\cal A}_{a}),\quad 
\Delta_{\varphi}\equiv \oplus\{\Delta_{\varphi}^{(a)}:a\in \sem\},\]
%
%
\begin{equation}
\label{eqn:counittwo}
\vep(x)\equiv 
\left\{
\begin{array}{cl}
0\quad &\mbox{ when }x\in \oplus \{{\cal A}_{a}:a\in \sem\setminus \{e\}\},\\
\\
\vep_{e}(x)\quad&\mbox{ when }x\in {\cal A}_{e}.
\end{array}
\right.
\end{equation}
Then $({\cal A}_{*},\delp,\vep)$ is a strictly proper counital C$^{*}$-bialgebra.
\end{Thm}

\noindent
We call $({\cal A}_{*},\Delta_{\varphi},\vep)$ in 
Theorem \ref{Thm:mainthree} by a (counital)
{\it C$^{*}$-bialgebra} associated with 
$\{({\cal A}_{a},\varphi_{a,b}):a,b\in \sem\}$.

We show a sufficient condition 
of the non-existence of the antipode for the C$^{*}$-bialgebra
in Theorem \ref{Thm:mainthree}.
For an element $a$ in a monoid $\sem$,
if there exists $b\in \sem$ such that $ba=e$,
then $a$ is called {\it left invertible}.
The multiplicative monoid $(\ndm{*},\boxtimes)$ in (\ref{eqn:nondegenerate})
has no left invertible element except the unit $1$.
%
%
\begin{lem}
\label{lem:noexist}
(Lemma 3.2 in \cite{TS02})
Assume that $\sem$ is a monoid such that
any element in $\sem\setminus\{e\}$ is not left invertible.
For a C$^{*}$-WCS $\{({\cal A}_{a},\varphi_{a,b}):a,b\in \sem\}$
over $\sem$ which satisfies (\ref{eqn:finiteness}),
let $({\cal A}_{*},\delp,\vep)$ be as in Theorem \ref{Thm:mainthree}
and let $(\tilde{{\cal A}}_{*},\hdelp,\tilde{\vep})$ be the smallest unitization
of $({\cal A}_{*},\delp,\vep)$ in Lemma \ref{lem:extensionunit}.
Then the antipode for any dense unital counital subbialgebra of 
$(\tilde{{\cal A}}_{*},\hdelp,\tilde{\vep})$ never exists.
\end{lem}

Let $\{({\cal A}_{a},\varphi_{a,b}):a,b\in \sem\}$ and
$\{({\cal B}_{a},\psi_{a,b}):a,b\in \sem\}$
be two C$^{*}$-WCSs over a monoid $\sem$.
Assume that $\{f_{a}:a\in \sem\}$ is a set of 
unital $*$-homomorphisms such that
$f_{a}\in {\rm Hom}({\cal A}_{a},{\cal B}_{a})$ for $a\in \sem$ and 
the following holds:
%
%
\begin{equation}
\label{eqn:morphism}
\psi_{a,b}\circ f_{ab}=(f_{a}\otimes f_{b})\circ \varphi_{a,b}
\quad(a,b\in \sem),\quad \vep_{{\cal B}_{e}}\circ f_{e}=\vep_{{\cal A}_{e}}.
\end{equation}
Then $f_{*}\equiv \oplus \{f_{a}:a\in \sem\}$
is a counital C$^{*}$-bialgebra morphism from ${\cal A}_{*}$ 
to ${\cal B}_{*}$.
If ${\cal B}_{a}={\cal A}_{a}$ and $f_{a}$ is bijective for each $a$,
then $f_{*}$ is a C$^{*}$-bialgebra automorphisms of ${\cal A}_{*}$.

Furthermore, the following lemma holds.
%
%
\begin{lem}
\label{lem:subalgebra}
(Lemma 3.3 in \cite{TS02})
For a C$^{*}$-WCS $\{({\cal A}_{a},\varphi_{a,b}):a,b\in \sem\}$,
if ${\cal B}_{a}$ is a unital $C^{*}$-subalgebra of ${\cal A}_{a}$
and the inclusion map $\iota_{a}$ of $B_{a}$ into ${\cal A}_{a}$
satisfies that
\[\varphi_{a,b}({\cal B}_{ab}) \subset {\cal B}_{a} \otimes {\cal B}_{b}
\quad(a,b\in \sem),
\quad 
\vep_{e}\circ \iota_{e}=\vep_{e},
\]
then ${\cal B}_{*}=\oplus\{{\cal B}_{a}:a\in \sem\}$ 
is a counital C$^{*}$-subbialgebra of $({\cal A}_{*},\delp,\vep)$. 
\end{lem}

Remark that Definition \ref{defi:axiom}
can be reformulated in the purely algebraic situation without involution.
For this ``algebraic weakly coassociative system",
Theorem \ref{Thm:mainthree} is also true as a bialgebra without unit.
Furthermore 
such a bialgebra is always extended to the unital bialgebra
by Lemma \ref{lem:extensionunittwo}.

%
%
\ssft{Proof of theorems}
\label{subsection:thirdfour}
We prove Theorem \ref{Thm:maintsfour} and \ref{Thm:maintwo} here.
\\

\noindent
{\it Proof of Theorem \ref{Thm:maintsfour}.}
By definition, we can verify that the set $\{\varphi_{A,B}:A,B\in \ndm{*}\}$ satisfies
the assumption in Theorem \ref{Thm:mainthree}
with respect to the monoid $\ndm{*}$ from the discussion 
in $\S$ \ref{section:second}.
Hence $\{(\coa,\varphi_{A,B}):A,B\in\ndm{*}\}$
is a C$^{*}$-weakly coassociative system.
From Theorem \ref{Thm:mainthree}, the statement holds.
\qedh

\noindent
{\it Proof of Theorem \ref{Thm:maintwo}.}
(i) From Lemma \ref{lem:extensionunit}, the statement holds.\\
(ii)
For $A\in \ndm{n}$,
let $s_{1}^{(A)},\ldots,s_{n}^{(A)}$ 
denote canonical generators of $\coa$ and 
let $\coa^{(0)}$ denote the dense $*$-subalgebra generated by them.
Let $\odot$ denote the algebraic tensor product.
By the definition of $\varphi_{A,B}$,
we see that $\varphi_{A,B}(\co{A\boxtimes B}^{(0)})\subset 
\co{A}^{(0)}\odot \co{B}^{(0)}$ for each $A,B\in\ndm{*}$.
Let ${\sf CK}^{(0)}_{*}$ denote the algebraic direct sum of 
the set $\{\coa^{(0)}:A\in\ndm{*}\}$.
If we write ${\goth A}_{0}$ as the unitization of ${\sf CK}^{(0)}_{*}$,
then the statement holds from Lemma \ref{lem:extensionunittwo}. \\
(iii)
From Lemma \ref{lem:noexist},
the statement holds.\\
(iv)
From $\S$ \ref{subsection:secondthree},
${\bf N}$ is a submonoid of $(\ndm{*},\boxtimes)$ such that any divisor 
of any element in ${\bf N}$
belongs to ${\bf N}$.
From this, $\delp({\sf C}_{*})\subset {\sf C}_{*}\otimes {\sf C}_{*}$.
Hence ${\sf C}_{*}$ is a C$^{*}$-subbialgebra of ${\sf CK}_{*}$.
By the definition of the comultiplication of $\co{*}$ in \cite{TS02},
${\sf C}_{*}$ is isomorphic to $\co{*}$ as a C$^{*}$-bialgebra.
\qedh

%
%
\sftt{Automorphisms and representations}
\label{section:fourth}
We show automorphisms and representations of 
$({\sf CK}_{*},\delp,\vep)$ in this section.
%
%
\ssft{Modified gauge action}
\label{subsection:fourthone}
For $A\in \ndm{n}$,
let $s_{1}^{(A)},\ldots,s_{n}^{(A)}$ denote canonical generators of $\coa$.
Define the action $\lambda^{(A)}$ of $U(1)$ on $\coa$ by
\[\lambda_{z}^{(A)}(s_{i}^{(A)})\equiv 
z^{\log n}s_{i}^{(A)}\quad(z\in U(1),\,i=1,\ldots,n).\]
Especially, $\lambda_{z}^{(1)}=id$ for each $z$.
We call $\lambda^{(A)}$ the {\it modified gauge action} of $\coa$.
Define the $*$-automorphism $\lambda^{(*)}_{z}$ of the C$^{*}$-algebra
${\sf CK}_{*}$ by
\[\lambda^{(*)}_{z}\equiv \oplus\{\lambda_{z}^{(A)}:
A\in\ndm{*}\}\quad(z\in U(1)).\]
Then we can verify that 
$\lambda^{(*)}_{z}$ is a C$^{*}$-bialgebra automorphism of $({\sf CK}_{*},\delp,\vep)$
for each $z$ and 
$\lambda^{(*)}_{z}\circ \lambda^{(*)}_{w}=\lambda^{(*)}_{zw}$
for each $z,w\in U(1)$.
Therefore $\lambda^{(*)}$ is an action of $U(1)$ on
the C$^{*}$-bialgebra $({\sf CK}_{*},\delp,\vep)$.
We call $\lambda^{(*)}$ the {\it gauge action} of ${\sf CK}_{*}$.

By definition,
the extension $\tilde{\lambda}^{(*)}$ of $\lambda^{(*)}$ 
on $\widetilde{{\sf CK}}_{*}$
is also an action of $U(1)$ on
the unital C$^{*}$-bialgebra $\widetilde{{\sf CK}}_{*}$.
Furthermore the restriction of $\tilde{\lambda}^{(*)}$
on ${\goth A}_{0}$ in Theorem \ref{Thm:maintwo} (ii)
is also an action of $U(1)$ on the bialgebra ${\goth A}_{0}$.
%
%
\ssft{Tensor product of permutative representations of 
Cuntz-Krieger algebras and its decomposition}
\label{subsection:fourthtwo}
In \cite{CKR01,CKR02,PFO01},
we introduced a class of representations of Cuntz-Krieger algebras.
We show that 
the tensor product associated with
the comultiplication $\delp$ of ${\sf CK}_{*}$
is closed on this class of representations.

Let $A,B,C\in \ndm{*}$ and 
let ${\rm Rep}\coa$ denote
the class of all unital $*$-representations of $\coa$.
For $(\pi_{1},\pi_{2})\in {\rm Rep}\coa\times {\rm Rep}\co{B}$,
define $\pi_{1}\ptimes \pi_{2}\in {\rm Rep}\co{A\boxtimes B}$ by
%
%
\begin{equation}
\label{eqn:definition}
\pi_{1}\ptimes \pi_{2}\equiv (\pi_{1}\otimes \pi_{2})\circ \varphi_{A,B}
\end{equation}
where $\varphi_{A,B}$ is as in (\ref{eqn:varphi}).
This defines the following operation:
\[\otimes_{\varphi}:{\rm Rep}\coa\times{\rm Rep}\co{B}
\to {\rm Rep}\co{A\boxtimes B}.\]
For $\pi_{1},\pi_{2}\in{\rm Rep}\coa$,
define the relation $\pi_{1}\sim \pi_{2}$ if $\pi_{1}$ and $\pi_{2}$ 
are unitarily equivalent.
For $\pi_{1},\pi_{1}^{'}\in {\rm Rep}\coa$,
$\pi_{2},\pi_{2}^{'}\in {\rm Rep}\co{B}$ and $\pi_{3}\in {\rm Rep}\co{C}$,
the following holds:
\begin{enumerate}
\item
If $\pi_{1}\sim \pi_{1}^{'}$ and $\pi_{2}\sim \pi_{2}^{'}$,
then $\pi_{1}\otimes_{\varphi} \pi_{2}\sim
\pi_{1}^{'}\otimes_{\varphi} \pi_{2}^{'}$.
\item
$\pi_{1}\otimes_{\varphi} (\pi_{2}\oplus \pi_{2}^{'})=
\pi_{1}\otimes_{\varphi} \pi_{2}\,\oplus\, \pi_{1}\otimes_{\varphi} \pi_{2}^{'}$.
\item
$\pi_{1}\otimes_{\varphi} (\pi_{2}\otimes_{\varphi} \pi_{3})
=(\pi_{1}\otimes_{\varphi} \pi_{2})\otimes_{\varphi} \pi_{3}$.
\end{enumerate}
Furthermore, by definition, we can verify that the following holds 
for $\ptimes$ in (\ref{eqn:definition})
and $\delp$ in Theorem \ref{Thm:maintsfour}:
\[(\pi_{1}\otimes \pi_{2})\circ \delp|_{\co{A\boxtimes B}}
=\pi_{1}\ptimes\pi_{2}
\quad ((\pi_{1},\pi_{2})\in{\rm Rep}\coa\times {\rm Rep}\co{B})\]
where $\co{A\boxtimes B},\coa$ and $\co{B}$
are naturally identified with C$^{*}$-subalgebras of ${\sf CK}_{*}$ and  
$\pi_{1}\otimes \pi_{2}$ is naturally identified
with a representation of ${\sf CK}_{*}\otimes {\sf CK}_{*}$.
From this,
the C$^{*}$-bialgebra structure of ${\sf CK}_{*}$
brings the new tensor product $\ptimes$. 
Since there exist representations $\pi_{1},\pi_{2}$ of ${\sf CK}_{*}$ such that
$\pi_{1}\ptimes \pi_{2}\not\sim \pi_{2}\ptimes \pi_{1}$,
$\ptimes$ is not symmetric.

Remark that the operation $\ptimes$ in (\ref{eqn:definition})
is defined for general representations of Cuntz-Krieger algebras.
We consider $\ptimes$ on the following small class of representations.
%
%
\begin{defi}
\label{defi:chara}
Let $A\in \ndm{n}$ and let $s_{1},\ldots,s_{n}$
denote the canonical generator of $\coa$.
A representation $({\cal H},\pi)$ of $\coa$
is permutative if there exists a complete orthonormal basis
$\{e_{l}\}_{l\in\Lambda}$ of ${\cal H}$,
a set $\{\Lambda_{i}\}_{i=1}^{n}$ of subsets of $\Lambda$
and a subset 
$\{m_{i,l}:i=1,\ldots,n,\,l\in \Lambda_{i}\}$ of $\Lambda$ such that
%
%
\begin{equation}
\label{eqn:first}
\pi(s_{i})e_{l}=
\left\{\begin{array}{ll}
e_{m_{i,l}}\quad &(l\in\Lambda_{i}),\\
\\
0\quad&(\mbox{otherwise}).\\
\end{array}
\right.
\end{equation}
\end{defi}
Roughly speaking,
a permutative representation preserves
a certain complete orthonormal basis.
For the operation $\ptimes$ in (\ref{eqn:definition}),
we see that 
if both $\pi_{1}$ and $\pi_{2}$ are permutative representations, 
then so is $\pi_{1}\ptimes \pi_{2}$.

According to \cite{CKR02},
any cyclic permutative representations 
of $\coa$ for $A\in \ndm{n}$ is characterized 
by a finite or an infinite word of $1,\ldots,n$.
Any permutative representation is decomposed into
cyclic permutative representations uniquely up to
unitary equivalence.
From this, the tensor product decomposition
of two cyclic permutative representations makes sense.
In consequence, the following statement holds.
%
%
\begin{prop}
\label{prop:deco}
For $A,B\in\ndm{*}$,
let $\pi_{1}$ and $\pi_{2}$ be cyclic permutative representations
of $\coa$ and $\co{B}$, respectively.
Then there exists a set $\{\eta_{n}:n\in\Omega\}$
of cyclic permutative representations of $\co{A\boxtimes B}$ such that 
%
%
\begin{equation}
\label{eqn:decomposition}
\pi_{1}\ptimes\pi_{2}=\bigoplus_{n\in\Omega}\eta_{n}
\end{equation}
and the cardinality of $\Omega$ is countably infinite at most.
Furthermore, this decomposition
is unique up to unitary equivalence.
\end{prop}
In particular,
we obtain the decomposition formulae in (\ref{eqn:decomposition})
for the case of Cuntz algebras \cite{TS01}.

%
%
\sftt{Examples of C$^{*}$-subbialgebras of ${\sf CK}_{*}$}
\label{section:fifth}
In this section,
we show examples of C$^{*}$-subbialgebras of ${\sf CK}_{*}$.
For $A\in \ndm{n}$,
let $s_{1}^{(A)},\ldots,s_{n}^{(A)}$ denote
canonical generators of $\coa$ in this section.
%
%
\begin{ex}
{\rm
Assume that a subset $\{A_{n}:n\geq 1\}$ of $\ndm{*}$
satisfies that $A_{n}\in \ndm{n}$ and 
\[A_{n}\boxtimes A_{m}=A_{nm}\quad(n,m\geq 1).\]
Then
\[\oplus\{\co{A_{n}}:n\geq 1\}\]
is a C$^{*}$-subbialgebra of ${\sf CK}_{*}$
if and only if ${\cal N}_{A_{n}}$ in (\ref{eqn:na})
is a subset of $\{(A_{m},A_{l}):m,l\geq 1\}$
for each $n\geq 1$.}
\end{ex}
%
%
\begin{ex}
{\rm
Fix $A\in \ndm{n}$ for $n\geq 2$.
Assume that if $A=B\boxtimes C$,
then $B$ or $C$ equals $1\in\ndm{1}$.
Then
\[\co{A^{*}}\equiv {\bf C}\oplus \co{A}
\oplus \co{A^{\boxtimes 2}}
\oplus \co{A^{\boxtimes 3}}\oplus \co{A^{\boxtimes 4}}\oplus \cdots\]
is a C$^{*}$-subbialgebra of ${\sf CK}_{*}$
where $A^{\boxtimes n}$ denotes
$A\boxtimes \cdots \boxtimes A$ ($n$-times).
}
\end{ex}
%
%
\begin{ex}
{\rm
For $A\in\ndm{n}$,
define the C$^{*}$-subalgebra $AF_{A}$ of $\coa$ by
\[AF_{A}\equiv C^{*}\langle\cup_{l\geq 1}
\{s_{J}^{(A)}(s_{K}^{(A)})^{*}:J,K\in\{1,\ldots,n\}^{l}\}\rangle\]
where $s_{J}^{(A)}\equiv s_{j_{1}}^{(A)}\cdots s_{j_{k}}^{(A)}$
when $J=(j_{1},\ldots,j_{k})$.
Define the C$^{*}$-subalgebra $AF_{*}$ of ${\sf CK}_{*}$ by
\[AF_{*}\equiv \oplus \{AF_{A}:A\in\ndm{*}\}.\]
Then $AF_{*}$ is a C$^{*}$-subbialgebra  of ${\sf CK}_{*}$
by Lemma \ref{lem:subalgebra}.
}
\end{ex}
%
%
\begin{ex}{\rm
In \cite{TS02},
we introduced various C$^{*}$-subbialgebras of $\co{*}$.
Since $\co{*}$ is a C$^{*}$-subbialgebra of ${\sf CK}_{*}$,
they are also C$^{*}$-subbialgebras of ${\sf CK}_{*}$.
}
\end{ex}
%
%
\begin{ex}
{\rm
According to $\S$ 6.4 in \cite{TS02},
we construct C$^{*}$-subbialgebras of ${\sf CK}_{*}$.
For $n\in {\bf N}$, $A\in \ndm{n}$ 
and a nonempty subset $\Sigma$ of $\{1,\ldots,n\}$,
define the C$^{*}$-subalgebra $\coa(\Sigma)$ of $\coa$ generated 
by the set $\{s_{i}^{(A)}:i\in \Sigma\}$ and the unit of $\coa$. 
For a set ${\bf \Sigma}\equiv \{\Sigma_{A}:A\in \ndm{*}\}$
such that $\Sigma_{A}\subset \{1,\ldots,n\}$
when $A\in \ndm{n}$,
define the C$^{*}$-subalgebra ${\sf CK}_{*}({\bf \Sigma})$ of ${\sf CK}_{*}$ by
\[{\sf CK}_{*}({\bf \Sigma})\equiv \oplus\{\coa(\Sigma_{A}):A\in \ndm{*}\}.\]
Define sets ${\bf 1}$, $\intercal$ and ${\bf 1}\cup \intercal$ by
${\bf 1}\equiv \{\Sigma_{A}:\Sigma_{A}=\{1\},\,A\in \ndm{*}\}$,
$\intercal\equiv \{\Sigma_{A}:\Sigma_{A}=\{n\}\mbox{ when }A\in\ndm{n}\}$
and ${\bf 1}\cup \intercal
\equiv \{\Sigma_{A}:\Sigma_{A}=\{1,n\}\mbox{ when }A\in\ndm{n}\}$.
Then we can verify that ${\sf CK}_{*}({\bf 1})$,
${\sf CK}_{*}(\intercal)$ and 
${\sf CK}_{*}({\bf 1}\cup \intercal)$ 
are counital non-commutative cocommutative C$^{*}$-subbialgebras 
of ${\sf CK}_{*}$ by Lemma \ref{lem:subalgebra}.
}
\end{ex}

%
%
\begin{ex}{\rm
We show a relation between 
the symbolic dynamical system and a commutative non-cocommutative 
C$^{*}$-subbialgebra of ${\sf CK}_{*}$.
For $A=(a_{ij})\in\ndm{n}$,
define the commutative C$^{*}$-subalgebra ${\cal C}_{A}$ of $\coa$ by
\[{\cal C}_{A}\equiv C^{*}\langle\{s_{J}^{(A)}(s_{J}^{(A)})^{*}:
J\in \nset{+}\}\rangle\]
where $\nset{+}=\bigcup_{l\in {\bf N}}\nset{l}$
and define the set 
\[X_{A}\equiv \{(j_{i})_{i=1}^{\infty}:a_{j_{i}j_{i+1}} = 1\}\] 
of all one-sided infinite sequences.
Then $X_{A}$ is a compact Hausdorff space with
respect to the relative topology of the product space
$\{1,\ldots,n\}^{\infty}$ of discrete topological spaces.
Let $C(X_{A})$ denote the C$^{*}$-algebra 
of all continuous complex-valued functions on $X_{A}$.
Then it is well-known that $C(X_{A})\cong {\cal C}_{A}$ \cite{CK}.
Especially, ${\cal C}_{1}\cong {\bf C}$.
Define the commutative C$^{*}$-subalgebra ${\cal SF}_{*}$ of ${\sf CK}_{*}$ by
%
%
\begin{equation}
\label{eqn:cantorone}
{\cal SF}_{*}\equiv \oplus\{{\cal C}_{A}:A\in\ndm{*}\}.
\end{equation}
Then ${\cal SF}_{*}$ is $*$-isomorphic to the C$^{*}$-algebra $C_{0}(X_{*})$
where $X_{*}$ denotes the direct sum of the set $\{X_{A}:A\in\ndm{*}\}$.
Furthermore,
we can verify that
${\cal SF}_{*}$ is counital 
commutative non-cocommutative C$^{*}$-subbialgebra of ${\sf CK}_{*}$.
We see that the smallest unitization $\widetilde{{\cal SF}}_{*}$ 
of ${\cal SF}_{*}$ is $*$-isomorphic to 
the C$^{*}$-algebra $C(\alpha X_{*})$
where $\alpha X_{*}$ denotes the one-point compactification of $X_{*}$.
}
\end{ex}


%
%

%
\end{document}